\documentclass{amsart}
\usepackage{amssymb}
\usepackage[]{amscd}
\usepackage{amsthm}
\usepackage{amstext}
\usepackage[all]{xy}
\newtheorem*{theo1}{Theorem 1}
\newtheorem*{theo2}{Theorem 2}
\newtheorem*{theo3}{Theorem 3}
\newtheorem{theo}{Theorem}
\newtheorem{prop}{Proposition}
\newtheorem{lem}{Lemma}[theo]

\newfont{\pr}{msbm10}
\newfont{\qr}{msbm8}
\newfont{\got}{eufm10}
\renewcommand{\P}{\hbox{\pr P}}
\newcommand{\I}{\mathcal I}
\newcommand{\p}{\hbox{\qr P}}
\renewcommand{\L}{\mathcal L}

\newcommand{\E}{\mathcal E}
\newcommand{\C}{\hbox{\qr C}}
\renewcommand{\O}{\mathcal O}
\newcommand{\opr}{ {\O}_{ {\p}^{r} } }
\newcommand{\res}[2]{{}_{#1}| #2} 

\begin{document}

\title {On the regularity of
varieties having an extremal secant line}

\author{Marie-Am\' elie Bertin}

\begin{abstract}
The Castelnuovo-Mumford regularity of
a non degenerate variety of degree $d$ and
dimension $n$ of ${\P}^{r}$ is conjectured to be
at most $d-r+n+1$ \cite{Cms}. This conjecture is
known to hold only in a few cases: curves
\cite{GLP}, smooth surfaces (\cite{P},\cite{L}),
and smooth varieties of codimension $2$ \cite{A}.
Varieties with an extremal secant line are
$(d-r+n)$-irregular and the classification of
smooth $(d-r+1)$-irregular curves in \cite{GLP}
shows that they all possess an extremal secant
line but the elliptic normal curve and some smooth rational curves
in ${\P}^{d-1}$. We give here a
complete classification of varieties having an 
extremal secant line and show that their regularity
is $d-r+n+1$. 
\end{abstract}

\email{mabertin@fourier.ujf-grenoble.fr}

\curraddr{Universit\' e de Grenoble I, Institut
Fourier, 38402 St Martin-d'H\` eres, FRANCE}
\maketitle
\section{Introduction}

Let $X$ be a complex $n$-dimensional variety, non degenerate and of degree
$d$ in ${\P}^{r}$.  
The $k$-regularity of $X$ has been defined by D.Mumford~\cite{Mu}
by  the following vanishing
$$H^{i}({\P}^{r},{\mathcal
I}_{X|{\p}^{r}}(k-i))=0\quad\quad\text{for all }
i\geq 1.$$
The $k$-regularity condition implies the $(k+1)$-regularity condition
, so that one defines the \textit{Castelnuovo-Mumford regularity} of $X$,
 as the the least integer $k$ such that $X$ is $k$-regular. We
denote by $reg(X)$ the regularity of $X$.
 This obscure vanishing condition find its origin in
G.Castelnuovo's work on the postulation of a variety $X$
\cite{C}; for the regularity of $X$ gives an upper bound on the 
least  integer $k$ for which hypersurfaces of degree $s\geq k$ cut on $X$
a complete regular linear system.
The regularity $X$ also encodes 
pieces of information on the syzygies of
the equations defining $X$.
Indeed, let $I_{X}$ be the saturated ideal of $X$
in $R=\hbox{\pr C}[x_{0},\ldots ,x_n]$ and
$$0\rightarrow L_{r+1}\rightarrow\cdots \rightarrow
L_{0}\rightarrow R\rightarrow R/I_{X}\rightarrow 0 $$
be a minimal free graded resolution of the $R$-module $R/I_{X}$;
the variety $X$ is $k$-regular if and only
if for all $i\geq 0$, one can find a base of $L_{i}$, which elements
are at most of degree $k+i$ \cite{BM}, definition 3.2.

This algebraic characterization of $k$-regularity
has an elementary geometric consequence:
the existence of a $k$-secant line to $X$ implies
that  $reg(X)\geq k$. 
One can deduce from Bertini's theorem on linear systems
that $X$ cannot have any $k$-secant line
for $k>d-r+n+1$; we will thus call \textit{extremal} 
a $(d-r+n+1)$-secant line to $X$. Varieties for which there
exist extremal secant lines are therefore $(d-r+n)$-irregular.
The \textit{regularity conjecture} (\cite{Cms},\cite{EG}) foresees that 
$reg(X)\leq d-r+n+1$. It has so far only been proved for curves
\cite{GLP}, smooth surfaces (\cite{P}, \cite{L}) and smooth $2$-codimensional
varieties \cite{A}.
Extremal varieties for the conjecture have only been classified for
curves by Gruson, Lazarsfeld and Peskine in \cite{GLP}. They are
smooth rational curves with an extremal secant line except for
the elliptic normal curve and some rational curves in ${\P}^{d-1}$.

Varieties having an extremal secant line thus form
a good sample for testing the regularity conjecture.

\textbf{Conventions and notations:}
We will use Hartshorne's book conventions on projective bundles.
Consider $n$ positive integers $1\leq a_1\leq \ldots \leq
a_n$, the smooth rational normal scroll $S(a_1,\ldots ,a_n)$ is
the tautological embedding of the projective bundle ${\P}(\oplus_{i=1}^{n}
{\O}_{{\p}^{1}}(a_i))$ in ${\P}^{n-1+\sum_{i=1}^{n}a_i}$; a non degenerate cone of vertex
${\P}^{k}$ over $S(a_1,\ldots a_n)$ in ${\P}^{n+k-\sum_{i=1}^{n}a_i}$  will be denoted by $S(\underbrace{0,\ldots
,0}_{k-\hbox{times}},a_1,\ldots ,a_n)$.  

Here we will prove the following results.

First of all we show the existence of
non degenerate $n$-folds $X$ of degree $d$ in ${\P}^r$,
having an extremal secant line, for any $d>r-n+1$, $n$ and $r\geq 2n+1$.

\begin{theo1} Let us fix an integer $n$, then choose $n$ non zero 
  positive integers $a_1\leq \ldots \leq a_n$ and a positive integer
$d> \sum_{i=1}^{n} a_i$.
Let us pick also two non zero positive integers $k_{1}$ and $k_{2}$ 
such  that $k_{1}+k_{2}=d-\sum_{i=1}^{n} a_i$ and an injection of
${\O}_{{\p}^{1}}$-module
$${\O}_{{\p}^{1}}(-k_{1})\oplus
{\O}_{{\p}^{1}}(-k_{2})\xrightarrow{\alpha}
{\O}_{{\p}^{1}}^{2}\oplus_{i=1}^{n}
{\O}_{{\p}^{1}}(a_{i}).$$
The cokernel map 
$${\O}_{{\p}^{1}}^{2}\oplus_{i=1}^{n}
{\O}_{{\p}^{1}}(a_{i})\xrightarrow{\beta} \oplus_{i=1}^{n}
{\O}_{{\p}^{1}}(b_{i})$$
gives a morphism of rational scrolls from $\overline{X}=S(b_1,\ldots ,b_n)$ to
$S=S(0,0,a_1,\ldots ,a_n)$ restriction of a linear projection
between their linear spans that induces an isomorphism of rational 
scrolls between $\overline{X}$ and its image $X$. The vertex of the 
cone $S$ is an extremal secant line for the $n$-fold scroll $X$.
Moreover any rational $n$-fold scroll having an extremal secant line
can be realized this way.
\end{theo1} 

This is completed by our classification theorem.

\begin{theo2} Let $X$ be a non degenerate variety
in ${\P}^{r}$ with $r-n\geq 2$ of degree $d$ and
dimension $n\geq 2$. Assume that $X$ has an extremal secant
line $l$. Then $X$ is either

\begin{itemize}
\item a cone over the Veronese
surface $V$ in ${\P}^{5}$ or

\item a cone over a smooth rational
scroll $\Sigma$, such that $l$ is  an
extremal secant line for $\Sigma$ or

\item a cone over the Veronese surface $V^{\prime}$ in 
${\p}^4$.
\end{itemize}
\end{theo2}
It is known that $reg(V)=2$ and $reg(V^{\prime})=3$.
From the algebraic characterization of $k$-regularity,
if $X$ is a cone over $\Sigma$ we have $reg(X)=reg(\Sigma )$,
so that to show that varieties with an extremal secant
line satisfy the regularity conjecture we only need to show it
for $n$-fold rational scrolls.

\begin{theo3} Let $X$ be a $n$-fold scroll over a
smooth curve $C$ of genus $g$, of degree $d$ in
${\P}^{r}$. Then we have the inequality
$$reg(X)\leq d-r+n+1.$$
\end{theo3}

\textbf{Acknowledgments} \quad This work was made while I was a Ph.D. student 
at Columbia University. Most of it constitutes my Ph.D. 
dissertation \cite{Ber} and grew out of discussions with my advisor Prof. Henry C. 
Pinkham, to whom I express my gratitude. The last theorem was inspired
by personal notes of Prof. Lazarsfeld \cite{L2},
kindly communicated to us by their author. 

\bigskip

\section{Rational scrolls
with an extremal secant line}

Beside the case where bisecant lines are
extremal secant lines for $X$, i.e. $X$ is
a variety of minimal degree, the existence of
varieties with extremal secant line  is not
clear. The latest construction of such varieties is due
to Lazarsfeld. He constructs in Remark 1 of
\cite{L} smooth surfaces of ${\P}^{r},r\geq 5$
with an extremal secant line as $1$-codimensional
subvarieties of a rational normal $3$-fold
scroll. The construction  readily
generalizes to higher dimension and furnishes
smooth varieties (actually rational scrolls) of
any degree $d>r-n+1$ in
${\P}^{r},r\geq 2n+1$ having an extremal secant
line. 

We explain here how to construct any rational
$n$-fold scroll $X$ that has an extremal secant
line. First we note that $X$ must be a $2$-codimensional
subvariety of a singular rational normal scroll.

\begin{prop} Let $X$ be a variety in ${\P}^{r}$
having an extremal secant line $l$. The image of
$X$ by the linear projection $\pi_{l}$ of center
$l$ is a $n$-dimensional variety $X^{\prime}$ of minimal
degree in ${\P}^{r-2}$. Moreover if
$X$ is a smooth rational scroll so is
$X^{\prime}$.
\end{prop}

\begin{proof}

If we project $X$ from this line $l$ to
${\P}^{r-2}$, we get a variety $X^{\prime}$ of
degree $d^{\prime}$ and dimension $n^{\prime}\leq
n$ in ${\P}^{r-2}$. Since $X$ is non degenerate,
$n=n^{\prime}$. It follows that $d^{\prime}=r-n-1$,
hence $X^{\prime}$ is of minimal degree.

The rational scroll $X$ is the image of
a smooth rational normal scroll $Z$ of degree $d$ and dimension $n$,
by a linear projection. 
We deduce from the classification of varieties of minimal
degree~(\cite{B},\cite{EH})
that $X^{\prime}$ is a rational normal scroll.
The induced rational map $g$ from $Z$ to
$X^{\prime}=\pi_l(X)$ is thus an elementary transformation (cf.\cite{M})
of $Z$ along $Y=\{ p_{1},\cdots ,p_{s}\}$.
It follows that $X^{\prime}$ is the tautological embedding of
a projective bundle, elementary transform of ${\P}(\E)$;
therefore $X^{\prime}$ is smooth.

\end{proof}
 
Thus our $n$-fold scroll lies on the
cone $<l,X^{\prime}>=S(0,0,a_{1},\cdots ,a_{n})$
for some integers $1\leq a_{1}\leq \cdots
\leq a_{n}$. We can reverse the process; on any
rational normal scroll $S=S(0,0,a_{1},\cdots
,a_{n})$, we can construct a smooth codimension
$2$ subscroll $X$ having the vertex of the cone
$S$ as extremal secant line.

\begin{theo} 
 Fix an integer $n$, then choose $n$ non zero 
  positive integers $a_1\leq \ldots \leq a_n$ and a positive integer
$d> \sum_{i=1}^{n} a_i$.
Let us pick also two non zero positive integers $k_{1}$ and $k_{2}$   
such  that $k_{1}+k_{2}=d-\sum_{i=1}^{n} a_i$ and an injection of
${\O}_{{\p}^{1}}$-module
$${\O}_{{\p}^{1}}(-k_{1})\oplus
{\O}_{{\p}^{1}}(-k_{2})\xrightarrow{\alpha}
{\O}_{{\p}^{1}}^{2}\oplus_{i=1}^{n}
{\O}_{{\p}^{1}}(a_{i}).$$
The cokernel map 
$${\O}_{{\p}^{1}}^{2}\oplus_{i=1}^{n}
{\O}_{{\p}^{1}}(a_{i})\xrightarrow{\beta} \oplus_{i=1}^{n}
{\O}_{{\p}^{1}}(b_{i})$$
gives a morphism of rational scrolls from  $\overline{X}=S(b_1,\ldots 
,b_n)$ to $S=S(0,0,a_1,\ldots ,a_n)$ restriction of a linear projection
between their linear spans that induces an isomorphism of rational 
scrolls between $\overline{X}$ and its image $X$. The vertex of the 
cone $S$ is an extremal secant line for the $n$-fold scroll $X$.
Moreover any rational $n$-fold scroll having an extremal secant line
can be realized this way.
\end{theo}
\begin{proof}
Let us denote by $\tilde S$ the projective bundle
${\P}({\O}_{{\p}^{1}}^2\oplus_{i=1}^{n}{\O}_{{\p}^{1}}(a_{i}))$ in ${\P}^{1}\times {\P}^r$,
where $r=n+1+\sum_{i=1}^{n} a_i$.
The second projection $\pi_2$ from the product ${\P}^{1}\times {\P}^r$ embeds
$\tilde S$ in ${\P}^r$ as a cone $S$ of vertex $l$ over the smooth 
rational normal scroll $X^{\prime}=S(a_1,\cdots ,a_n)$. Its restriction to $\tilde S$ 
is a desingularization  of $S$ of exceptional locus $E={\P}^1\times l={\P}({\O}_{{\p}^1}^2)$.
A simple Chern class computation shows that $\sum b_i=d+n-1$, with $b_i\geq 1$ for all $1\leq i\leq n$.
Let us denote by $\overline{X}$ the smooth rational normal scroll 
$S(b_1,\cdots ,b_n)$ in ${\P}^{d+n-1}$.
Since $h^0({\P}^1, {\O}_{{\p}^{1}}(-k_1)\oplus{\O}_{{\p}^{1}}(-k_2))=0$, the map $\alpha$
defines a morphism $f$ from $\overline{X}$ to $S$, which image $X$ is a non degenerated rational 
subscroll of $S$ in ${\P}^r$. Moreover since $\overline{X}$ is linearly normal 
and $S$ non degenerated, $f$ is the restriction to $\overline{X}$ of a linear projection. 
Since the restriction $\beta_2$ of $\beta$ to 
$\sum_{i=1}^{n}{\O}_{{\p}^{1}}(a_{i})$ has to be an injection of 
${\O}_{{\p}^1}$-modules, it defines an elementary transformation between
$S(a_1,\cdots ,a_n)$ and $\overline{X}$, induced by a linear projection
from ${\P}^{d+n-1}$ to ${\P}^{r-2}$, which factors through $f$ and the projection of ${\P}^r$ 
from the line $l$. Therefore $X$ is smooth of dimension $n$, so that 
$f$ induces an isomorphism of rational scrolls between $\overline{X}$ and 
$X$ and $deg X=d$.

We denote by $\tilde X$ the strict transform of $X$ in $\tilde
S$. Let $\sigma$ denote the embedding map ${\P}^{1}\rightarrow l\subset {\P}^{r}$. 
Notice that the diagonal $\Delta=\{ (t,p)\in E :p=\sigma(t)\} $ in $E$ is isomophic to $l$ by $\pi_2$,
so that the schemes $Z=\tilde X\cap\Delta$ and $X\cap l$ are isomorphic by $\pi_2$. 
For each $i=1,2$, the injections induced by $\alpha$,
${\O}_{{\p}^1}(-k_1)\xrightarrow{\alpha_{1,i}} {\O}_{{\p}^1}^2 $, 
define in $E$ over ${\P}^1$ a rational curve $C_i={\P}({\O}_{{\p}^1}(k_i))$. 
Let $\oplus \hbox{\pr C}^{n_i} w_i$ denote the cokernel of
$${\O}_{{\p}^{1}}(-k_1)\oplus{\O}_{{\p}^{1}}(-k_2)
\xrightarrow{\alpha_1=\alpha_{1,1}+\alpha_{1,2}}{\O}_{{\p}^{1}}^2;$$
Remark that $ 1\leq n_i\leq 2$.  
Since for each $i=1,\ldots ,s$, the line $
\pi_1^{-1}(w_i)\cap E$ meets $\Delta$, the scheme $Z$ is
supported on $\{ (w_1,\sigma(w_1)),\cdots ,(w_s,\sigma (w_s))\}$ in $\Delta$.
If $n_i=2$, the smoothness of each curve $C_i$ at each point of
$C_i\cap \Delta$, shows that $(w_i,\sigma (w_i))$ is on $C_i\cap \Delta$ for
both curve $C_1$ and $C_2$. The cotangent space of $\tilde X\cap E$ at 
$(w_i,sigma (w_i))$ is then the direct sum of the cotangent spaces of $C_i$ at
this point. Therefore the multiplicity of $\tilde X\cap \Delta$ at
$(w_i,\sigma (w_i))$ is $2$. From this, we deduce that the length of $Z$ is
$\sum n_i=d-r+n+1$. 
\end{proof}

\section{The classification of varieties having
an extremal secant line}

Let $X$ be a $n$-dimensional variety in
${\P}^{r}$ having an extremal secant line $l$.
From Gruson, Lazarsfeld and Peskine's
classification of $(d-r+1)$-irregular curves \cite{GLP}, if
$n=1,r\geq 3$, $X$ is a \emph{smooth rational
curve}. This and Del~Pezzo-Bertini's theorem
\cite{B},\cite{EH} that classifies varieties of minimal degree, are
the key ingredients for our classification theorem. 

\begin{theo}
Let $X$ be a non degenerate variety
in ${\P}^{r}$ with $r-n\geq 2$ of degree $d$ and
dimension $n\geq 2$. Assume that $X$ has an extremal secant
line $l$. Then $X$ is either
\begin{itemize}
\item a cone over the Veronese
surface $V$ in ${\P}^{5}$ or

\item a cone over a smooth rational
scroll $\Sigma$, such that $l$ lies  is  an
extremal secant line for $\Sigma$ or

\item a cone over the Veronese surface $V^{\prime}$ in 
${\p}^4$.
\end{itemize}
\end{theo}
\begin{proof}
Let $\tilde X\xrightarrow{\phi} X$ be a desingularization of $X$
of exceptional locus $E$.
We want to determine the linearly normal variety
$\overline{X}$ of which $X$ is projection. We show by induction on
$n$ that  $\overline{X}$ is a variety of minimal
degree $d$. To do so, we need the existence of an hyperplane 
section of $X$ containing $l\cap X=\{ p_1,\ldots ,p_s\}$ that is desingularisezd by $\phi$.

\vfill
\eject

\begin{lem} \hfill

\begin{enumerate}
\item The variety $X$ is smooth  in a
neighborhood of $X\cap l$.

\item The variety $\tilde X$ satisfies
$h^{1}(\tilde X,\O_{\tilde X})=0$.
\end{enumerate}
\end{lem}

\begin{proof}[Proof of (1)]

This is true for $n=1$ since then $X$ is a
smooth rational curve \cite{GLP}.
Let us assume it is true for $k<n$.

Consider the linear system $\L$ of hyperplanes $H$
containing $l$. The generic element of the restriction to $X$ of $\L$
satisfies  our induction hypothesis so that any
$p_{i}$ in $X\cap l$ is  a smooth point of $h=X\cap H$ for $H$
generic. Since any smooth point $p$ on a Cartier divisor $D$ of a variety
$X$ is a smooth point of $X$, $p_i$ is a smooth point of $X$.
\end{proof}

\begin{proof}[Proof of (2)]
Since
$h^{1}(\tilde X,\O_{\tilde X})$ is a
birational invariant for smooth varieties,
we deduce that
$$h^{1}(\tilde X,\O_{\tilde X})=h^{1}(\tilde
X^{\prime},\O_{\tilde X^{\prime}})$$
for any desingularization $\tilde X^{\prime}$ of
the variety $X^{\prime}$. Moreover since
$X^{\prime}$ is a variety of minimal degree,
hence a cone over a smooth rational variety, we have
$h^{1}(\tilde X^{\prime},\O_{\tilde
X^{\prime}})=0$.
\end{proof}

Let us now prove that $X$ is the image by a
regular projection of a variety of minimal degree
$\overline{X}$.

\begin{lem} Let $H$ be a generic
hyperplane section of
$X$. Recall that $\phi\colon\tilde
X\rightarrow X$ is a desingularization of $X$.
\begin{enumerate}

\item A generic element of the linear system
$|\phi^{*}(\O_{X}(H))|$ is smooth and the
dimension of this system is
$d+n-1$.

\item The rational map $f: \tilde
X\dashrightarrow {\P}^{d+n-1}$ defined by this
system is in fact regular and maps $\tilde X$
onto a variety $\overline{X}$ of dimension
$n$ and degree
$d$.
\end{enumerate}
\end{lem}
\begin{proof}
For $n=1$, \cite{GLP}, the curve $X$ is a
smooth rational curve of degree $d$, hence is a
regular projection of some rational normal curve
in ${\P}^{d}$.

For $1\leq k<n$ assume that for any 
k-dimensional varieties $X$ having an extremal
secant line, the total transform  $|\O_{\tilde
X}(H)|$ of the linear system of hyperplane
sections of $X$, gives a map from
$\tilde X$ to ${\P}^{d+k-1}$ which image $\overline{X}$ is
a $k$-dimensional variety of degree $d$, that
projects onto $X$.

Let us consider on $\tilde X$ the total transform 
$\L$ of the system cut out on $X$ by the system
of hyperplanes containing $l$. Let $D$ be a generic hyperplane section 
of $X$ containing $X\cap l$; a
generic element of $D$ is irreducible and $\phi^{*}D$ of $\L$ is smooth.

Indeed $\phi^{*}D$ is smooth away from its
base locus by Bertini's theorem  on
singularities of a generic member of a linear system. 
The variety $X$ is smooth
at $p_{1},\cdots p_{s}$, hence the base locus
of $\L$ consists of the $s$ points 
$\phi^{-1}p_{1},\cdots ,\phi^{-1}p_{s}$.
The divisor $D$ is an irreducible hypersurface of
$X$ by  Bertini's irreducibility theorem ($r-2\geq 2$);
it has an extremal secant line by construction,
hence $D$ is  smooth in a neighborhood of
$p_{1},\cdots ,p_{s}$ by  lemma 2.1 (1).
The Cartier divisor $\phi^{*}D$ is then smooth
at $\phi^{-1}p_{1},\cdots ,\phi^{-1}p_{s}$.
This also implies that a generic element of
$|\O_{\tilde  X}(\phi^{*}H)|$ is smooth.

We can thus find an irreductible and non degenerated 
hyperplane section of $X$, $H_{0}$
such that $\phi^{*}H_{0}$ is a smooth element of
$\L$.
We have the following exact sequence:
$$0\rightarrow \O_{\tilde X}\rightarrow
\O_{\tilde X}(\phi^{*}H)\rightarrow
\O_{\phi^{*}
H_{0}}((\phi^{*}H)|_{\phi^{*} H_{0}})
\rightarrow 0.$$

By lemma 2.1 (2), 
$h^{0}(\O_{\tilde X}(\phi^{*}H))=1+h^{0}(\O_{\phi^{*}
H_{0}}((\phi^{*}H)|_{\phi^{*} H_{0}}) )$.  The restriction of $\phi$ 
to the closure of $\phi^{*}H_{0}\setminus E$ in $\tilde X$ is a desingularization of 
$H_{0}$ so that we can apply the induction hypothesis to a generic
hyperplane section $h$ of $H_{0}$ to get:
$$h^{0}(\tilde H_{0},
\O_{\tilde H_{O}}(\phi|_{\tilde H_{0}}^{*}h))
=d+n-1.$$
Since a generic element of $|\O_{\phi^{*}
H_{0}}(\phi^{*} H|_{\phi^{*} H_{0}})|$ doesn't
meet $E$,
$$h^{0}(\phi^{*}
H_{0},\O_{\phi^{*} H_{0}}(\phi^{*}
H|_{\phi^{*} H_{0}}))=h^{0}(\tilde H_{0},
\O_{\tilde H_{0}}(\phi^{*} H|_{\tilde H_{0}}))
=d+n-1.$$

Since $|\O_{\tilde X}(\phi^{*}H)|$ is base point free, the degree of $\overline{X}$ is 
$d$ and the rational map $\psi$ that $|\O_{\tilde X}(\phi^{*}H)|$ defines is
regular.

The dimension  of $\overline{X}$ is at
least $n-1$. Since $X$, is the union of the points in the pencil
formed by $2$ independant hyperplane sections $H_1$ and $H_2$ of $X$ containing $l\cap 
X$, $\overline{X}$ is cut out by the pencil generated by the
images of the strict transforms of $H_1$ and $H_2$ in $\tilde X$ and therefore has dimension $n$.

The linear system $|\O_{\tilde 
X}(\phi^{*}H)|$ generically separates points, so that the map $\psi$ is  
birational onto its image.
 Its inverse $\psi^{-1}$ is then a rational map from $\overline{X}$
onto $X$ that extends to a linear projection $\pi$ from
${\P}^{d+n-1}$ onto ${\P}^{r}$, since $\overline{X}$
is linearly normal. The map $\psi^{-1}$ is regular  since 
$\overline{X}$ and $X$ have the same dimension and the same degree.
\end{proof}

We can now use the Del Pezzo-Bertini theorem to
conclude. Assume that $r-n=2$, then $\overline{X}$ is a cone over the Veronese
surface $V$ in ${\P}^5$ and $X$ is the linear projection of 
$\overline{X}$ by a point. Since the extremal secant line $l$ meets $X$ at 
smooth points, it must be the image of a $3$-secant $2$-plane of 
$V$, hence the center of projection lies in the linear span of $V$
and $X$ is a cone over the Veronese surface in ${\P}^4$.  
If $r\geq n+3$ the variety $\overline{X}$ can be neither a cone
over the Veronese surface nor the Veronese
surface itself, unless $\overline{X}=X$.
Indeed, if it were so, $X$ would be a cone
$<L^{k},V^{\prime}>$ over the generic projection
of the Veronese to ${\P}^{4}$ or  over
the Steiner surface in ${\P}^{3}$, with
$k=n-3$. So clearly, we would have $r\leq n+2$.
If $n=2$, we cannot have any projection
of the Veronese surface either.

Hence $\overline{X}$ must be a rational scroll unless $X$
was already a variety of minimal degree. Since
the projection $\pi$ is regular, if
$\overline{X}$ is smooth, $X$ is a smooth
rational scroll. If $\overline{X}$ is a cone
$<L^{k},S^{n-k-1}>$ over a smooth scroll  of
minimal degree $S^{n-k-1}$, $X$ must be the cone
$<L^{k},\pi (S^{n-k-1})>$, hence is a cone over a smooth rational
scroll. Any extremal secant line moreover has to
lie in the linear space generated by $S^{n-k-1}$,
for it has to meet $X$ at smooth points.
\end{proof}

\section{The regularity of smooth scrolls}

In this section we  prove that varieties having an extremal secant 
line are $d-r+n+1$-regular. By the classification theorem, we only need to prove
the regularity of smooth scrolls having an extremal secant line.

This will follow from this more general result on the regularity
of smooth scrolls over a smooth curve $C$.

\begin{theo}
Let $X$ be a smooth $n$-dimensional scroll over a smooth curve $C$ 
of genus $g$,
embedded in ${\P}^{r}$ as a non degenerate variety of degree $d$.
The regularity of $X$ is bounded by $d-r+n+1$, as predicted
by the regularity conjecture.
\end{theo}

The proof is closely related to Gruson, Lazarsfeld and Peskine's proof
of the regularity conjecture for curves. As they do, we first describe
$X$ as a projection from a product of a variety $\mathcal S$ that is scheme
theoretically a vanishing locus by a Beilinson type construction~\cite{OSS}; in our case
$\mathcal S$ is the the ruled subvariety of $C\times {\P}^{r}$ which projects 
onto the $n$-fold scroll $X$.
We get a set-theoretical 
description of $X$ as a degeneracy locus, in the same manner as in \cite{GLP}.
We can no longer conclude as in \cite{GLP}, for the ideal sheaf of $X$ and of
the degeneracy locus no longer differ along a 
$0$-dimensional scheme, but along $d-r+n$ fibers of $X$. We will
see that this degeneracy locus is fibered over $C$ and obtained from $X$ by "replacing" 
$d-r+n$ fibers of $X$ by multiple fibers and this extra
information will be enough to bypass this problem.
\begin{proof}
The $k$-regularity of $X$ is equivalent to the $k-1$-normality of $X$ 
and the  vanishing of $H^i({\O}_{X}(k-1-i))$. 
Note that by degeneration of the Leray spectral sequence
$$H^{i}(X,{\O}_{X}(d-r+n-i))\simeq H^{i}(C,\pi_{*}({\O}_{X}(d-r+n-i))),$$
so that we only have to show that 
$$H^{i}({\P}^{r},{\I}_{X|{\p}^{r}}(d-r+n+1-i))=0\,\,\,\hbox{\rm for } i=1,2.$$ 

The embedding ${\P}(E)\xrightarrow{p} X\subset {\P}^r={\P}(V)$ 
factors through an embedding $h$ of ${\P}(E)$ in $C\times {\P}^r$ 
and the second projection $f$ onto ${\P}^r$, hence
corresponds to the data of an exact sequence of ${\O}_{C}$-modules
$$ 0\rightarrow M \xrightarrow{i_p} V\otimes_{\C} {\O}_{C}\rightarrow E 
\rightarrow 0$$
such that $H^1(C,M)=0.$
Let us denote by $\pi$ the first projection from the product 
$C\times {\P}^r$.
From the Euler exact sequence on ${\P}^{r}$, twisted by 
${\O}_{{\p}^r}(1)$ and pulled back by $f$ to $C\times {\P}^r$ and from the pull back by $\pi$ of 
the morphism  $i_p$ we get a morphism $\pi^{*}M \xrightarrow{s} f^{*}(\opr (1))$ of ${\O}_{C\times 
{\p}^r}$-modules which vanishing locus is scheme-theoretically 
$h({\P}(E))$.

The Koszul complex associated to $\pi^{*}(M)\otimes f^{*}(\opr 
(-1))\xrightarrow{s\otimes f^{*}(\opr (-1))} {\O}_{C\times {\p}^{r}}$
gives a locally free resolution of
$h_{*}({\O}_{{\p}(E)})$:

$$
\begin{matrix}
0\rightarrow\pi^{*}(\wedge^{r-n+1} M)\otimes f^{*}(\opr (n-r-1)) \rightarrow
\cdots \rightarrow \pi^{*}(\wedge^{j} M)\otimes f^{*}(\opr (-j)) \rightarrow
\cdots \\
\rightarrow \pi^{*}(\wedge^{2}M)\otimes f^{*}(\opr (-2)) \rightarrow 
\pi^{*}(M)\otimes f^{*}(\opr (-1)) \rightarrow \opr \rightarrow h_{*}({\O}_{\p 
(E)})\rightarrow 0\\
\end{matrix}.$$

Indeed this complex is exact since $s\otimes f^{*}(\opr (-1))$ is generically surjective
and its cokernel $h_{*}({\O}_{{\p}(E)})$ is of expected codimension 
($rk(M)-1=r-n$) in the locally Cohen-Macaulay scheme ${\O}_{C\times 
{\p}^{r}}$ (\cite{Sz} p.~344).

Since the article of Gruson, Lazarsfeld and Peskine, it is now a standard trick ,
to twist this resolution by the pull back by $\pi$ of a 
suitable line bundle $A$ on $C$, in order to be able to push down to 
$X$ in ${\P}^r$ the pieces of information of 
this resolution of $h({\P}(E))$. 

\begin{lem}
There is a line bundle $A$ on $C$ which satisfies
$$h^{1}(M\otimes A)=0, \,\,\,\,\, h^{1}(\wedge^{2}M\otimes A)=0, \,\,\,\,\, 
h^{1}(A)=0.$$ For such a line bundle we have then $h^{0}(M\otimes 
A)=(r-n)(d-r+n)+1$ and $h^{0}(A)=d-r+n$.
\end{lem}
\begin{proof}
  
The proof of this is similar to \cite{GLP} (lemma 1.7). 
First we have a strictly decreasing filtration 
by vector bundles
$F_{i}$ such that $\frac {F_{i}}{F_{i+1}}=L_{i}.$
is a line bundle of negative degree
$$ M=F_{1}\supset F_{2}\supset \cdots \supset F_{r+1-n} \supset 0=F_{r+2-n}$$
In order to get $H^{1}(C,M\otimes A)=0$ and $H^{1}(C,\wedge^{2}M\otimes A)=0$, it is enough to choose $A$ 
such that $H^{1}(C,L_{i}\otimes A)=0$ and $H^{1}(C,L_{i}\otimes L_{j}\otimes A)=0$.

Since a generic line bundle of degree $\geq g-1$ is non special,
it is enough to take $A$ generic, so that $deg(A) +r-n-1-d\geq g-1$,
that is to say $deg(A)\geq g-r+n+d=d-codim(h({\O}_{{\p}(E)}))+g\geq g-1$.
Indeed we have $deg(L_{i}\otimes {L_j})\geq -d+r-n-1$.

Since $rk(M)=r+1-n$  the smallest possible degree
for $deg(L_{i})$ is $r-n-d$, so that we also have $H^{1}(C,M\otimes A)=0$,.
Note that  $h^{1}(C,A)=0$, so that by Riemann-Roch's 
theorem $h^{0}(C,A)=d+n-r+1$. Since $h^{1}(C,M\otimes A)=0$ and $h^{1}(C,L_{i}\otimes A)=0$ 
we have $h^{0}(C,M\otimes A)=(r-n)(d+n-r)+1$.
\end{proof}
Applying K\" unneth's formula to the push forward by $f$ of the Koszul resolution of $h_{*}({\O}_{{\p}(E)})$ 
twisted  by $\pi^{*}(A)$ we get the following complex
\begin{multline*}
0\rightarrow H^0(\wedge^{r-n+1}M\otimes A)\otimes_{\C}\opr 
(n-r-1)\rightarrow \cdots \rightarrow H^{0}(\wedge^{2}M\otimes 
A)\otimes_{\C }\opr (-2) \\
\rightarrow H^{0}(M\otimes A )\otimes_{\C} \opr (-1)\xrightarrow{u} 
H^{0}(A)\otimes_{\C }\opr \rightarrow
{p}_{*}({h}^{*}{\O}_{{\p}(E)}\otimes \pi^{*}(A))\rightarrow 0. 
\end{multline*}
Let $t$ denote the $C$-projective bundle structure map of ${\P}(E)$.
For the choice of line bundle $A$ as in the previous 
lemma, by the same argument as in (\cite{GLP}, (1.5)), this complex is exact.
We get $f_{*}\circ h_{*}(t^{*}A)\simeq 
f_{*}(h_{*}({\O}_{{\p}(E)})\otimes \pi^{*}A)$ by the projection formula.

The degeneracy locus of the morphism $\opr^{(r-n)(d+n-r)+1} (-1)\xrightarrow{u} 
\opr^{d+n-r+1}$ that we obtain this way is set-theoretically  $X=p(h({\O}_{{\p}(E)}))$.
The cokernel of $u$ is isomorphic to $p_{*}t^{*}A$ by  the projection formula applied to 
$p=f\circ h$ and $t=\pi\circ h$. 

A generic section $\sigma$ of $A$ induces an exact sequence
$$ 0 \rightarrow{\O}_{C} \xrightarrow{\sigma}  A \rightarrow A\otimes 
{\O}_{D} \rightarrow 0$$
whith $D=\mathop{\sum}\limits _{i=1}^{d+n-r} p_{i}$ and where the $p_i$'s are distinct points of $C$.
The push forward by $p$ of this exact sequence is
$$\begin{CD}
0 @>>> p_{*}{\O}_{{\p}(E)} @>>> p_{*}t^{*}A @>>> {\mathcal F} @>>> R^{1}p_{*}
{\O}_{{\p}(E)}=0,
\end{CD}$$
where $\mathcal F $ denotes $p_{*}(t^{*}(A\otimes {\O}_{D}))$. 

Moreover since $p$ is a birational
morphism and $X$ is smooth, $p_{*}({\O}_{{\p}(E)})\simeq {\O}_{X}$.
On each fiber $t^{-1}(p_{i})$, we have
$t^{*}(A)_{p_{i}}\simeq {\O}_{{\p}(E),p_{i}}$. Therefore we get
$t^{*}(A\otimes {\O}_{D})\simeq t^{*}{\O}_{D}$, so 
that ${\mathcal F} \simeq \mathop{\oplus }\limits_{i=1}^{d-r+n}{\O}_{X_{i}}$ 
where $X_{i}$ is the fiber of the scroll $X$ over the point $p_{i}$.

As in \cite{GLP} (proof of theorem 2.1 p.~500) the existence of the commutative diagram 
$$
\begin{CD}
&&&&&& 0&& 0 &\\
&&&&&&@VVV @VVV &\\
&&0 &&&&\langle t^{*}\sigma \rangle \otimes {\O}_{{\p}^{r}} @>>> 
{\O}_{X}
@>>> 0\\
& & @VVV && @VVV @VVp_{*}t^{*}\sigma V &\\
0@>>> K @>>> {\O}_{{\p}^{r}}^{l}(-1) @>u>>
H^{0}(t^{*}A)\otimes_{\C} {\O}_{{\p}^{r}}
@>>> p_{*}t^{*}A @>>> 0
\\ && @VVV @| @VVV @VVV &\\
0@>>> N @>>>{\O}_{{\p}^{r}}^{l}(-1) @>v>>
{\O}_{{\p}^{r}}^{d+n-r} @>>> {\mathcal F} @>>> 0\\
&& @VVV && @VVV @VVV &\\
  &&{\I}_{X|{\p}^r}&&&& 0 &&0 &&\\
&& @VVV && &&  &&\\
&&0&&&&&&\\
\end{CD},$$
which allows us to compare the regularity of $X$ with
the regularity of $N$, follows from the minimality of 
$$0\rightarrow K=ker(u)\rightarrow {\O}_{{\p}^{r}}^{l}(-1) \xrightarrow{u}
H^{0}(t^{*}A)\otimes_{\C} {\O}_{{\p}^{r}}
\rightarrow p_{*}t^{*}A. $$
In \cite{GLP} the morphism corresponding to $u$ is constructed from a  
minimal free resolution; the local minimality we seek here is taken 
care of by the following lemma.

\begin{lem}(A minimality criterion)
Let $L_{\cdot }$ be the Koszul complex associated to some map
$L_{1} \xrightarrow{s} R$ where $R$ is a local ring of maximal ideal 
$\hbox{\got n}$.
Let us assume that there is a local morphism
$S\hookrightarrow R$ turning $R$ into a finite free $S$-module of 
rank $l$.
Any $R$-module $N$ has via this morphism an $S$-module structure that 
we denote by $\res{S}{N}$.
Consider the map of $S$-modules $\res {S}{L_{1}} \xrightarrow{\phi} S$,
obtained by composing $\res {S}{L_{1}}\xrightarrow{\res{S}{s}} R$ with 
$R \xrightarrow{b\wedge -} S\simeq \wedge^{l}R$, where $b$ is the
wedge product of the first $h-1$ elements of a base of the free $S$-module $\res{S}{L_{1}}$.
The Koszul complex $K_{\cdot }$ of $S$-modules associated to $\phi$
surjects onto the complex of free modules $\res {S}{L_{\cdot }}$.
In particular this last complex is minimal if $K_{\cdot }$ is minimal.
Moreover if $\res {S}{s}(L_{1})\subset \hbox{\got m}R$, where 
$\hbox{\got m}$ is the
maximal ideal of $S$, then the complex $L_{\cdot }$ is minimal.
\end{lem}
\begin{proof}
The natural surjection of $R$-modules
$$\begin{CD}
R\otimes_{S}\res {S}{L_{1}} @>\psi>> L_{1}\\
x\otimes l&\longmapsto &xl\\
\end{CD}$$
gives a commutative diagram
$$\begin{CD}
L_{1}@>s>> R\\
@A\psi AA @AA\wr A\\
R\otimes_{S}\res {S}{L_{1}} @>R\otimes_{S}\phi >> R\\
\end{CD}$$
that induces a surjective map of complexes of free $R$-modules
between $K_{\cdot }\otimes_{S}R$ and $L_{\cdot}$.
Indeed by construction the map of $S$-modules $b\wedge -$ is 
surjective hence $(b\wedge -)\otimes_{S}R$ is also surjective. It
thus induces an automorphism $\sigma$ of $R$. It is then clear that $s\circ \psi =
{\sigma }^{-1}\circ R\otimes_{S}\phi$. The surjection between the complexes of $R$-modules induced by 
$\psi$ gives a surjection $K_{\cdot}\rightarrow \res{S}{L_{\cdot }}$ as claimed.

If $\res{S}{s}(\res{S}{L_{1}})\subset \hbox{\got m}R $ after 
wedging by $b$, we have $\phi(\res {S}{L_{1}}
)\subset \hbox{\got m}S$ so that $\res{S}{ L_{1}} \xrightarrow{\phi} S$
is minimal when $\res{S}{L_{1}} \xrightarrow{\res{S}{s}}\res{S}{R}$ 
is minimal.
\end{proof}
By K\" unneth formula $f_{*}(\pi^{*}(A))=H^{0}(A)\otimes_{\C}\opr$,
so that $f_{*}\pi^{*}A$ is a finite locally free $\opr$-module of rank 
$h^{0}(A)$. We are exactly in the situation of the criterion since locally, the complex $K_{\cdot }(s)\otimes \pi^{*}A$ 
is simply the Koszul complex associated to $\pi^{*}(M)_{x}\otimes 
{\pi}^{*}A_{x} \rightarrow \pi^{*}A_{x}$. Since the map 
$H^{0}(A)\otimes_{\C} {\O}_{{\p}^{r}} \rightarrow p_{*}t^{*} A$
induces a surjection on global sections, none of the sections
defining $H^{0}(A)\otimes_{\C} {\O}_{{\p}^{r}} \rightarrow p_{*}t^{*}A$ vanishes.
This shows then that $0\rightarrow K\rightarrow {\O}_{{\p}^{r}}^{l}(-1) \xrightarrow{u}
H^{0}(t^{*}A)\otimes_{\C} {\O}_{{\p}^{r}}\rightarrow p_{*}t^{*}A$ is 
locally minimal.

\begin{lem}
We have $H^{2}(K(k))=0$ for all $k>-n$.
\end{lem}
\begin{proof}
This is equivalent to show that $H^{1}(C\times {\p}^{r},{\I}_{{\p}(E)|{\p}^{r}}\otimes\pi^{*}A\otimes 
f^{*}\opr (k))=0$ for all $k>-n.$
Indeed since $f$ is projective and only has fibers of dimension $1$, 
$R^{j}f_{*}$ vanishes on any coherent sheaf for $j\geq 2$, so that
$$R^{1}f_{*}(\pi^{*}(M\otimes A)\otimes f^{*}(\opr (k-1)))\rightarrow  
R^{1}f_{*}(\I_{{\p}(E)|{\p}^{r}}\otimes\pi^{*}A\otimes f^{*}\opr 
(k))$$ is surjective.
By non speciality of $A$ and use of the K\" unneth formula we get 
$R^{1}f_{*}(\I_{{\p}(E)|{\p}^{r}}\otimes\pi^{*}A\otimes 
f^{*}\opr (k))=0$.
\end{proof}
 
Consider now the Koszul resolution of $h_{*}({\O}_{{\p}(E)}$
twisted by $f^{*}({\O}_{{\p}^{r}}(k))$
\begin{multline*}
\cdots \rightarrow \pi^{*}(\wedge^{2}M\otimes A)\otimes f^{*}(\opr(k-2)) 
\rightarrow \pi^{*}(M\otimes A)\otimes f^{*}(\opr (k-1))\rightarrow \\
\pi^{*}A\otimes f^{*}(\opr (k)) \rightarrow h_{*}{\O}_{{\p}(E)}\otimes 
f^{*}(\opr (k))\rightarrow 0.
\end{multline*}

Let us denote by ${\mathcal F}_{j}$ the cokernel of
$$\begin{CD} 
\pi^{*} (\wedge^{j+1} M\otimes A)\otimes f^{*}(\opr (k-j-1))@>>> 
\pi^{*}(\wedge^{j}M\otimes A)\otimes f^{*}(\opr (k-j)).
\end{CD}$$
\begin{lem}
We have 
$R^{i}\pi_{*}(\I_{{\p}(E)|{\p}^{r}}\otimes\pi^{*}A\otimes f^{*}\opr 
(k))=0$ for $i\geq 1$.
\end{lem}
\begin{proof}
Since the fibers of $\pi$ are of dimension $r$ and $\pi$ is
projective we know that $R^{r+1}\pi_{*}$ vanishes on any coherent 
sheaf, so that $R^{r+1}\pi_{*}$ is right exact.
Then
$R^{r}\pi_{*}(\pi^{*}(\wedge^{j} M\otimes A)\otimes 
f^{*}\opr(k-j))=\wedge^{j}M\otimes A\otimes H^{r}(\opr ,\opr 
(k-j))=0$ for all $ k>-n$, since $j\leq rk(M)=r+1-n$ and $k-j> -r-1$.
Therefore $R^{r}{\mathcal F}_{j}=0$ for all $j$. Repeating this argument, we 
deduce
that $R^{i}\pi_{*} {\mathcal F}_{j}=0$ for all $j$ when $k> -n$ and
$H^{1}(C,\pi_{*}(\pi^{*}(M\otimes A)\otimes f^{*}(\opr 
(k)))=H^{1}(C,M\otimes A)\otimes H^{0}(\opr ,\opr (k))=0.$
Moreover $H^{2}(C,\pi_{*}{\mathcal F}_{1})=0$, since 
$\pi_{*}{\mathcal F}_{1}$ is 
coherent ($\pi$ is projective) and $dim C=1$.
Therefore $H^{1}(\pi_{*}{\mathcal F}_{0})=H^{1}(\pi_{*}(\I_{{\p}(E)|{\p}^{r}}\otimes\pi^{* 
}A\otimes f^{*}\opr (k)))=0$.
\end{proof}

To show that $X$ is $d-r+n+1$-regular, it suffices to show that $H^{i}(N(k))=0$ 
for $k=d-r+n+1-i$ and $i=1,2$.

\begin{lem}  We have $H^{i}(N(k))=0$, for $k=d-r+n+1-i$ and $i=1,2$.
\end{lem}

We deduce from the local minimality of 
$$0\rightarrow K=ker(u)\rightarrow {\O}_{{\p}^{r}}^{l}(-1) \xrightarrow{u}
H^{0}(t^{*}A)\otimes_{\C} {\O}_{{\p}^{r}}
\rightarrow p_{*}t^{*}A $$
that the following complex extracted from the previous commutative diagram 
$$\begin{CD}
0 @>>> N @>>>{\O}_{{\p}^{r}}^{l}(-1) @>v>> {\O}_{{\p}^{r}}^{d+n-r} 
@>>> 0\end{CD} \qquad (*)$$ 
is also locally minimal.

Now the sheafification of a minimal resolution
of $\bigoplus_{l\in {\mathbb Z}} H^{0}({\mathcal F}(l))$ is of the form
$$
  0\rightarrow {\O}^{n_{r-1}}_{{\p}^{r}}(-r-1)\rightarrow \cdots 
\rightarrow 
\opr^{(r+1-n)(d+n-r)}(-1)\xrightarrow{w}
\opr^{d-r+n} \rightarrow {\mathcal F} \rightarrow
0.$$

Since the exact sequence $(*)$ comes from the sheafification
of a locally free minimal resolution of $\oplus_{l\geq 
l_{0}}H^{0}({\mathcal F }
(l))$ for some $l_{0}\geq 0$, we can construct the following 
commutative diagram
$$\begin{CD}
&&0&&0&&&&\\
&&@VVV @VVV &&&&\\
0@>>> N @>>> \opr^{l}(-1) @>v>> Im (v) @>>> 0\\
&& @VVV @VVV @VV\wr V &&\\
0@>>> P @>>> \opr^{(r+1-n)(d-r+n)}(-1) @>w>> Im w @>>> 0\\
&&@V\phi VV @VVV &&&&\\
&&\opr^{d+n-r-1} (-1) @>\simeq >> \opr^{d+n-r-1} (-1)&&&&\\
&&@VVV @VVV &&&&\\
&&0&&0&&&&\\
\end{CD} $$
Since $P$ has a resolution of type
$$\begin{CD}
0@>>>{\O}^{n_{r+1}}_{{\p}^{r}}(-r-1) @>>> \cdots 
@>>>{\O}^{n_{2}}_{{\p}^{r}}(-2)@>w_1>> P@>>>0\end{CD},$$
we have $H^{i}({\P}^{r},P(d-r+n+1-i))=0$ for $i=1,2$ ( by lemma 1.6 of
\cite{GLP}).

Therefore to show that $H^{i}({\P}^{r},N(d-r+n+1-i))=0$ for $i=1,2$, it
suffices to show that $ker (q)$ is $d-r+n+1$-regular, where $ker(q)$ 
is defined by the following trianglular diagram
$$\xymatrix{
{\opr^{n_2}(-2) } \ar[r]^{w_{1}} \ar[dr]_{q} & P \ar[d]^{\phi }\ar[r] 
&0\\
&\opr^{d+n-r-1}(-1) \ar[d] & \\
& 0 & }$$

This holds by the now standard trick due to Gruson, Lazarsfeld and 
Peskine, consisting in applying Lemma 1.6 of 
\cite{GLP} to the Eagon-Northcott complex associated to $q$ twisted by 
${\O}_{{\p}^r}(d+n-r-1)$.
\end{proof}


\begin{thebibliography}{19}
\bibitem{A} 
A. Alzati, \emph{A new
Castelnuovo bound for two codimensional
subvarieties in ${\P}^{r}$}, Proc. Amer. J.
Soc., {\textbf 114}, (1992), p. 607-611.

\bibitem{BM} D. Bayer and D. Mumford, \emph{
What can be computed in  algebraic geometry?},
Computational algebraic geometry and
commutative algebra (Cortona 1991), p. 1-48,
Sympos. Math. XXXIV, Cambridge Univ.
Press, Cambridge (1993).


\bibitem{Be} A. Beauville, Surfaces alg\'
ebriques complexes, Ast\' erisque, vol. 54 ,
Soci\' et\' e  Math\' ematique de France, (1978).

\bibitem{Ber} M.~A. Bertin, \textit{On the regularity of varieties having
an extremal secant line}, Ph.D. thesis, Columbia University, (2000).

\bibitem{B}
 E. Bertini, Introduzione alla
geometria proiettiva degli iperspazi, Enrico
Spoeri, Pisa , (1907).

 

\bibitem{C} G. Castelnuovo, \textit{Sui
multipli di una serie lineare di gruppi di punti
appartenente ad una curva algebrica}, rend.
Circ. Mat. Palermo, {\bf 7} (1893).

\bibitem{Cms} \bysame , Memorie Scelte, Nicola Zanichelli, Bologna 1937



\bibitem{EG} D.~Eisenbud and S.~Goto,
\emph{Linear free resolutions and  minimal
multiplicity}, J.Alg., \textbf{ 88}, (1984) p.
89-133

\bibitem{EH} D.Eisenbud and J.Harris, \emph{
Varieties of minimal degree (a centenial
account)},  Bowdoin conference in Algebraic
Geometry (Brunswick, Maine 1985), Proceedings of
Symposia in Pure Mathematics XXXXVI,
 part 1, p. 1-13, Am. Math. Soc., Providence RI (1987).


\bibitem{GLP} L.~Gruson, R.~Lazarsfeld,
Ch.~Peskine, \emph{ On a theorem of Castelnuovo,
and the equations defining space curves}, Inv.
Math , {\textbf 72}, (1983)






\bibitem{L} R.~Lazarsfeld, \emph{ A
Castelnuovo bound for smooth surfaces}, Duke J.
1987

\bibitem{L2} R.~Lazarsfeld, \emph{Private
notes, June 1982, communicated in May 1999}.


\bibitem{M} M.~Maruyama, \emph{Elementary
transformations in the theory of algebraic
vector bundles}, Algebraic Geometry (La Rabida
1981), (1982) p. 241-266 

\bibitem{Mu} D. Mumford, Lectures on curves
on an algebraic  surface, Annals of Math.
studies, vol.  59,  Princeton Univ. Press, (1966)


\bibitem{OSS} C.Okonek, M.Schneider and
H.Spindler, Vector Bundles on Complex Projective
Spaces,  Progress in Math , vol. 3 , Birkha\"
user, (1980) 

\bibitem{P} H.C.Pinkham, \emph{Castelnuovo
bound for smooth surfaces}, Invent.Math., (1986)
p. 321-332




\bibitem{Sz} L.~Szpiro, \emph{Sur les travaux
de Kempf, Kleinman et Laskov sur les diviseurs
exceptionnels},  S\' eminaire Bourbaki, expos\' e
no. 417, Lecture Notes in Math, {\textbf 317},
Springer V.

\end{thebibliography}
\end{document}